\documentclass[12pt]{amsart}
 \numberwithin{equation}{section}
 
\theoremstyle{plain}
\newtheorem{thm}{Theorem}[section]
\newtheorem{lemma}[thm]{Lemma}
\newtheorem{pro}[thm]{Proposition}

\newtheorem{ex}[thm]{Example}
\newtheorem{de}[thm]{Definition}

\begin{document}

\title[RUELLE OPERATOR FOR INFINITE CONFORMAL IFS]%
{RUELLE OPERATOR FOR INFINITE CONFORMAL IFS}
 \author[Xiao-Peng Chen, Li-Yan Wu and Yuan-Ling Ye ]{Xiao-Peng Chen, Li-Yan Wu and Yuan-Ling Ye*}
\address{Xiao-Peng Chen:  School of Mathematics and Statistics, Huazhong University of Science and Technology and
School of Mathematical Sciences, South China Normal University,
Wuhan 430074, PR China} \email{chenxiao002214336@yahoo.cn}
\address{Li-Yan Wu: School of Mathematical Sciences, South China Normal University and Teaching and Research Section of Mathematics
Public and Basic education department Guangdong Vocational College
of Industry and Commerce, Guangzhou 510510 PR China}
\email{gzwuliyan@tom.com}
\address{Yuan-Ling Ye: School of Mathematical Sciences, South China Normal University, Guangzhou 510631, PR China}
 \email{ylye@scnu.edu.cn}

\keywords{Conformal map, iterated function system, open set
condition, Ruelle operator, equicontinuity.}
\thanks{{\it 2000 Mathematics Subject Classification}: 28A78, 28A80, 47B38.}
 \date{\today}
\thanks{* The corresponding author: the Research was partially supported by SRF for ROCS, SEM }


\maketitle

\bigskip

\begin{center}
{\bf Abstract}
\end{center}

{\small Let  $(X, \{ w_j \}_{j=1}^m, \{ p_j \}_{j=1}^m)$ ($2 \leq m
< \infty$) be a contractive iterated function system (IFS), where
$X$ is  a compact subset of ${\Bbb{R}}^d$. It is well known that
there exists a unique nonempty compact set $K$ such that
$K=\bigcup_{j=1}^m w_j(K)$. Moreover, the Ruelle operator  on $C(K)$
determined by the IFS  $(X, \{ w_j \}_{j=1}^m, \{ p_j \}_{j=1}^m)$
($2 \leq m < \infty$) has been introduced in \cite{FL}.  In the
present paper, the Ruelle operators determined by the infinite
conformal IFSs are discussed. Some separation properties for the
infinite conformal IFSs are investigated  by using the Ruelle
operator.}

\bigskip
\bigskip



\bigskip

\section{\bf{Introduction}}

\bigskip
\bigskip

 The (finite) conformal iterated function systems (IFSs) have been studied
in various ways. Fan and Lau \cite{FL} considered the eigenmeasure
by using the Ruelle operator. Fan et al. (\cite{FL}, \cite{LRY})
extended Schief's idea \cite{S} to prove the equivalence of the OSC
and SOSC.  Peres et al. \cite{PRSS} showed that for conformal IFSs
satisfying the H\"{o}lder condition with the invariant set $K$, both
 OSC and SOSC are equivalent to $0< H^s( K) < \infty$, where $s$
is the zero point of the pressure function $P(\cdot)$. A simple
proof of the result was also given by Ye \cite{Ye}.

The infinite  iterated function systems (IFSs) of similarity maps
were developed by Moran \cite{M}. Moran \cite{M} extended the
classical results for finite systems of similitudes satisfying the
open set condition to the infinite case. Mauldin and Urba\'{n}ski
\cite{MU} considered the dimension and Hausdorff measure of the
limit sets of the infinite conformal IFSs. They found a way to
determine the Hausdorff dimension of the limit sets under the
condition that the  basic open sets were connected. In this paper,
we continue the study of  the infinite conformal IFSs. We only
assume that the space $X$ is locally connected. This is similar to
the finite case\cite{Ye}.

Let $(X, \{w_j\}_{j=1}^m, \{p_j\}_{j=1}^m)$ be  an  IFS, where $
w_j$'s are contractive self-maps on a compact subset $X \subset
{\Bbb{R}}^d$ and  $p_j$'s are positive Dini functions on $X$
\cite{FL}. Then there exists a unique compact set  $K$ invariant
under  IFS, i.e.,
$$K=\bigcup_{j=1}^m w_j(K).$$
The Ruelle operator on $C(K)$ is  defined by
\begin{eqnarray} \label{0.1}
T(f)(x)=\sum\limits_{j=1}^mp_j(w_j(x)) f(w_j(x)),\quad f \in C(K),
\end{eqnarray}
where $C(K)$ is the space of all continuous functions on $X$. For
the importance of studying the Ruelle operators can be found in
\cite{LY}. We could generalize the Ruelle operators to the infinite
conformal IFSs. Then we obtain the  PF-property(Theorem \ref{0.3})
using the Ruelle operator we define.

The separation properties(OSC and SOSC) are useful for studying the
(finite) conformal IFSs. However, for the infinite conformal IFSs,
Moran \cite{M} showed that self-similar set generated by a countable
system of similitudes may not be s-set even if the {\it open set
condition} (OSC) is satisfies. Szarek and Wedrychowicz \cite{SW}
 showed that for infinite IFSs, the OSC does not imply the SOSC. So it is necessary to consider the weak
 separation properties for the  infinite
conformal IFSs. Let $\{s_j \}_{j=1}^\infty$ be an infinite conformal
IFS on an open set $U_0$. We call the infinite conformal IFS $\{s_j
\}_{j=1}^\infty$ satisfies the {\it finite open set condition} if
for any integer $n$, there is a nonempty open set $U_n \subset U_0$
such that $s_i(U_n)\subset U_n$ for any  $ i \leq n$ and
$s_i(U_n)\cap s_j(U_n)=\emptyset$ for any $ i, j\leq n,i\neq j$.
Such weak separation property has been used  to study the Hausdorff
dimension of invariant set\cite{M}.
 We give a sufficient condition for   the infinite conformal
IFS $\{s_j \}_{j=1}^\infty$ satisfies the {\it finite (strong) open
condition} in Section \ref{0.2}.  Our basic results are:

\begin{thm} \label{0.3}
 Let $(X, \{s_j\}_{j=1}^\infty, \{ p_j \}_{j=1}^\infty)$ be an infinite conformal uniformly Dini IFS,
and let $J$ be the limit set of the  IFS $\{s_j \}_{j=1}^\infty$.
Then there exists  a unique $0<h \in C(\overline J)$ and a unique
probability measure $\mu \in M(\overline J)$ such that
\begin{eqnarray*}Th=\varrho h, \quad T^*\mu =\varrho \mu, \quad <\mu,
h>=1.
\end{eqnarray*}
Moreover, for every $f\in C(\overline J)$, $\varrho^{-n}T^nf$ converges
to $<\mu, f>h$ in the supremum norm, and for every $\xi \in
M(\overline J)$, $\varrho^{-n}T^{*n}\xi$ converges weakly to $<\xi,
h>\mu$.
\end{thm}

\bigskip

\begin{thm} \label{0.7} Let $\{s_j \}_{j=1}^\infty$ be an infinite conformal IFS, and let $J$ be
 the limit set of the IFS.
 If Hausdorff measure $ H^a (J)=H^a ( \overline J)>0$,
then the system $\{s_j \}_{j=1}^\infty$ satisfies the finite strong
open set condition.
\end{thm}

\bigskip

 Theorem \ref{0.3} generalizes THEOREM 1.1 of the paper
\cite{FL}.  Theorem \ref{0.7} extends the result of the paper
\cite{Ye}. Observe that the system is infinite, to prove Theorem
\ref{0.3} we need to assume that the potential functions are
uniformaly Dini continuous. Under the assumption that there exists a
number $a$ determined by the Ruelle operator, we  can prove Theorem
\ref{0.7} by making use of Schief's idea \cite{S}.

The paper is organized as follows. In Section 2, we present some
concepts for infinite conformal  IFSs and some elementary facts
about the Ruelle operators determined by infinite conformal IFSs. In
Section 3, we study the Ruelle operators defined by the infinite
conformal IFSs.
 In Section 4, we investigate the separation properties of the infinite conformal
 IFSs.


\bigskip

\section
{\bf Preliminaries}

\bigskip
\bigskip

In the paper, we always assume that $X$ is a locally connected compact subset of ${\Bbb{R}}^d$. Let $w$ be a self-map on compact set $X \subset {\Bbb{R}}^d$. We call $w$ a {\it conformal map} if there exists an open set $X \subset U_0$ such that $w$ is continuous differentiable on $U_0$ and for each $x \in U_0$, $w'(x)$ is a self-similar matrix.
We call the $\{s_j \}_{j=1}^\infty$ is an {\it infinite conformal iterated
function system} (IFS), if each $s_j (1 \leq j <\infty)$ is self-conformal map.
The family of functions $\{ p_j \}_{j=1}^\infty$ is said to be
{\it uniformly Dini continuous} on $U_0$ provided that $\int_0^1\frac {\alpha (t)} t <+\infty$, where
$$\alpha (t):=\sup\limits_{j \geq 1}\sup\limits_{ x, y \in U_0 \atop |x-y| \leq t}| \log
p_j(x)-\log p_j(y)|$$ and $\sum \limits_{j=1}^\infty p_j(x)< \infty$
on $U_0$. Under this assumption, we call the triple $(X,
\{s_i\}_{i=1}^\infty, $ $\{ p_j \}_{j=1}^\infty)$ an {\it infinite
conformal uniformaly Dini IFS}

Throughout the paper, we always assume that the infinite conformal IFS $\{ s_j \}_{j=1}^\infty$ satisfies the following conditions:

(i) each $s_j: X \rightarrow X$ is one-to-one and the system $\{s_j
\}_{j=1}^\infty$ is uniformly contractive, i.e., there exists
$0<s<1$ such that for any $1 \leq j <\infty$, $|s_j(x)-s_j(y)| \leq
s|x-y|, \quad \forall x,y \in X$;

 (ii) each ${s_j}$ is conformal on $U_0 (\supset X)$ with $0<\inf\limits_{x \in
U_0}|s'_j(x)| \leq \sup\limits_{x \in U_0} |s'_j(x)|<1$;

 (iii)
$\{|s'_j(\cdot)|\}_{j=1}^\infty$ is uniformly Dini continuous on $U_0$.

\bigskip

Let the $\{s_j \}_{j=1}^\infty$ be an infinite conformal IFS. Denote
$I=\{1, 2, 3, \cdots \}$,
\[ \begin{split} &I^n=\{{\mathbf u}=u_1 u_2\cdots u_n: u_j \in I, \ \forall 1 \leq j \leq n
\},\\
&I^\infty=\{{\mathbf v}=v_1 v_2\cdots : v_j \in I,  \ \forall 1 \leq
j < \infty \},\end{split} \] and let $I^*=\bigcup\limits_{n \geq 1}
I^n$. For arbitrary ${\mathbf w}=w_1 w_2 \cdots w_k \in I^k$, we let
$|{\mathbf w}|=k$ and $s_{\mathbf {w}} = s_{w_1} \circ s_{w_2}
\cdots \circ s_{w_k}$.

For any ${\mathbf w} \in {I^* \bigcup I^\infty}$ and integer $n$, if
$n \leq |{\mathbf w}|$, we denote by ${\mathbf w} |_n $ the word
$w_1 w_2 \cdots w_n$. By the contractiveness of the system $\{ s_j
\}_{j=1}^\infty$, we have
$$\lim_{n \rightarrow \infty} diam\{ s_{{\mathbf w}|_n}(X) \} \rightarrow 0,\quad \forall {\mathbf
w} \in {I^\infty}.$$ This implies the set $\bigcap
\limits_{n=1}^{\infty}s_{{\mathbf w}|_n}(X)$ is a singleton. Hence
we can define a map $\pi: I^\infty \rightarrow X$ by
$$\pi({\mathbf w}): =\bigcap
\limits_{n=1}^{\infty}s_{{\mathbf w}|_n}(X).$$
We endow the $I^\infty$ with the metric $d({\mathbf w}, {\mathbf u})=e^{-n({\mathbf w}, {\mathbf
u})}$, where $n({\mathbf w}, {\mathbf u})$ is the largest integer
$n$ such that ${\mathbf w}|_n={\mathbf u}|_n$. We can prove that the
map $\pi$ is continuous.

Define
\begin{equation} \label{1}
J=\pi(I^\infty)=\bigcup
\limits_{{\mathbf w} \in I^\infty} \bigcap \limits_{n=1}^\infty
s_{{\mathbf w} |_n}(X).
\end{equation}
We can check that the set $J$ satisfies the following equation:
     $$J=\bigcup _{j \in I}s_j(J).$$
The set $J$ is said to be the limit set of the system \cite{MU}. For
 an infinite conformal uniformaly Dini IFS $(X, \{s_i\}_{i=1}^\infty, $ $\{ p_j
 \}_{j=1}^\infty)$,
 we define an operator $T : C(\overline J) \to C(\overline J)$ by
\begin{equation} \label{2}
Tf(x) = \sum_{j=1}^\infty p_j (s_j(x)) f(s_j (x)).
\end{equation}
$T$ is called the Ruelle operator of the system.  The dual operator
$T^{\ast}$ on the measure space $M(\overline J)$ is given by
\begin{eqnarray} \label{3}
T^{\ast} \mu(E)=\sum_{j=1}^\infty \int_{s_j^{-1}(E)} p_j(x) d \mu(x)
\quad \mbox{ for any Borel set } E \subseteq \overline J
\end{eqnarray}
 (For the finite see e.g. \cite{BDEG}).

For ${\mathbf j} = j_{1}j_{2} \cdots j_{n},$  we let
$$p_{s_{\mathbf
{j}}}(x) = p_{j_1}(s_{j_2} \circ s_{j_3} \circ \cdots \circ
s_{j_{n}}(x)) \cdots p_{j_{n-1}}(s_{j_n}(x)) p_{j_{n}}(x).$$ By
induction we have for any integer $n$,
\begin{eqnarray*}
T^{n}f(x)=\sum\limits_{|{\mathbf j}|=n}p_{s_{\mathbf
j}}(x)f(s_{\mathbf j}(x)).
\end{eqnarray*}
Let $\varrho = \varrho(T)$ be the spectral radius of $T$.  Since $T$ is a
positive operator, we have $\|T^{n}1\|=\| T^{n} \|$  and
\begin{eqnarray*}
\varrho = \lim\limits_{n\rightarrow\infty} \| T^{n} \|^{
\frac{1}{n}}=\lim\limits_{n\rightarrow\infty} \| T^{n}1 \|^{
\frac{1}{n}}.
\end{eqnarray*}

For the sake of convenience, we
 can assume that
 $$diam (U_0)=\sup\{|x-y|:x, y\in U_0 \} \leq 1.$$

Let $J$ be the limit set of the system $\{ s_j \}_{j=1}^\infty$. For
any $r>0$, denote $B(J,r)=\bigcup \limits_{x\in J}B(x,r)$. Choose
$\kappa>0$ such that
\begin{equation} \label{eq2.1}
X_0 :=\bigcup \limits_{x\in X}\overline{ B(x, \kappa)} \subset U_0.
\end{equation}
 From the contractiveness of $s_j$'s, we can
show that $s_j(X_0) \subset X_0$ for any $j \in I$. Hence we may
assume that
\begin{eqnarray*}
B(J,\kappa)\subset X.
\end{eqnarray*}

\bigskip

 We call the infinite conformal IFS
$\{s_j \}_{j=1}^\infty$ satisfies the {\it open set condition} (OSC) if
there exists a nonempty bounded open set $U \subset U_0$ such that
$$s_j(U)\subset U  \ \mbox{ and } \
s_i(U) \cap s_j(U)=\emptyset \quad \mbox{ for any } i, j \in I
\mbox{ and } i \neq j.$$ Such a open set $U$ is called a {\it basic
open set} for the system $\{s_j \}_{j=1}^\infty$. If moreover $U
\bigcap J \neq \emptyset$, then the system $\{ s_j \}_{j=1}^\infty$
is said to satisfy the {\it strong open set condition} (SOSC).

To study the separation properties of the infinite conformal IFS, we
introduce the {\em finite open set condition} as follows.

\bigskip

\begin{de} \label{2.1.} We call the infinite conformal IFS
$\{s_j \}_{j=1}^\infty$ satisfies the finite open set condition if
for any integer $n$, there is a nonempty open set $U_n \subseteq
U_0$ such that $s_j(U_n) \subset U_n \ \mbox{ for each } 1 \leq j
\leq n$ and $s_i(U_n) \cap s_j(U_n)=\emptyset \ \mbox{ for each } 1
\leq i, j \leq n \mbox{ and } i \neq j$. The finite strong open set
condition holds if furthermore $U_n \bigcap J \neq \emptyset$.
\end{de}

\bigskip

Obviously the finite open set condition is weaker than the OSC.

\bigskip



\bigskip

\section
{\bf RUELLLE OPERATOR}

\bigskip

 \begin{pro} \label{3.1}  Let $\{ s_j \}_{j=1}^\infty$ be an infinite conformal IFS,
 and let $J$ be the limit set of the system $\{ s_j \}_{j=1}^\infty$. Then
     $$\overline {\{s_{\mathbf {u}}(x): {\mathbf u} \in I^*\}} = \overline J, \quad \forall x\in \overline
     J.$$
 \end{pro}

\bigskip

 \noindent{\bf Proof.} For any $x \in \overline J$
and ${\mathbf u} \in I^*$, we have
 $$s_{\mathbf u}(x) \in s_{\mathbf {u}}(\overline J) \subset \overline{ s_{\mathbf {u}} (J)} \subset \overline
 J.$$
 This implies that $\overline {\{s_{\mathbf {u}}(x): {\mathbf u} \in I^*\}} \subset \overline J$.

 For any $y \in J$, by  noticing that (\ref{1}) and the contractiveness of the system $\{ s_j \}_{j=1}^\infty$,
 there exists  ${\mathbf w} \in I^*$ such that
 $$\lim_{n \to \infty}\{s_{\mathbf w|_n}(x)\}=y, \quad \forall x \in \overline J.$$
 This implies that $J \subset \overline {\{s_{\mathbf {u}}(x): {\mathbf u} \in I^*\}}$.
 We have
$$\overline J \subset \overline {\{s_{\mathbf {u}}(x): {\mathbf u} \in I^*\}}.$$
Thus we conclude the assertion.
   \hfill$\Box$

   \bigskip

\begin{pro}\label{P3.2} Let the operator $T$ be defined as in
(\ref{2}). Then

{ \em (i)} $ \min_{x \in \overline J} \varrho^{-n}T^{n}1(x) \leq 1
\leq \max_{x \in \overline J} \varrho^{-n}T^{n}1(x), \quad \forall n
> 0$.

{\em (ii)} if there exist $\lambda > 0$ and $0 < h \in C(\overline J)$ such that
$T h=\lambda h$, then $\lambda = \varrho$ and there exist $ A, B > 0$ such that
$$
A \leq \varrho^{-n}T^{n}1(x) \leq B, \quad \forall  n > 0.
$$
\end{pro}

\medskip

\noindent{\bf Proof. } (i) We will prove the second inequality of
(i), the first inequality  is similar. Suppose it is not true, then
there exists a $k$ such that
 $\|T^k 1\| < \varrho^k$. Hence
$$
\varrho  =\big(\varrho(T^k)\big)^{\frac{1}k} \leq \| T^k \|^{\frac{1}{k}} = \|T^k
1\|^{\frac 1k} <\varrho,
$$
which it is a contradiction.

(ii) Let $a_{1} = \min_{x \in \overline J}
h(x)$, $a_{2}=\max_{x \in \overline J}h(x).$  Then
\begin{equation} \label{P3.2-1}
0< \frac {a_{1}} {a_{2}} \leq \frac{h(x)}{a_{2}}=
\frac{\lambda^{-n}} {a_{2}}T^{n}h(x) \leq \lambda^{-n}T^{n}1(x) .
\end{equation}
Since $T^n$ is a positive operator, we have $$\frac {a_{1}} {a_{2}}
\leq \lambda^{-n} \|T^{n}\|.$$
Similarly we can show that
\begin{equation} \label{P3.2-2}
\lambda^{-n}T^{n}1(x) \leq \frac {a_2} {a_1}.
\end{equation}
This implies that
$$\lambda^{-n} \| T^{n} \| \leq \frac {a_{2}} {a_{1}}.$$
Hence $\varrho =\lim\limits_{n \rightarrow
\infty}\|T^{n}\|^{\frac{1}{n}}=\lambda$. This, together with
(\ref{P3.2-1}) and (\ref{P3.2-2}), implies that
$$A:= \frac {a_{1}} {a_{2}} \leq \varrho^{-n}T^{n}1(x) \leq B:=\frac {a_2} {a_1}, \quad \forall  n > 0. $$
\hfill$\Box$
\bigskip

We call the operator $T: C(\overline J) \to C(\overline J)$  {\it irreducible} if for any
non-trivial, non-negative $f \in C(\overline J)$ and for any $x \in \overline J$,  there exists $n>0$ such that $T^{n}f(x)>0$.

\bigskip

\begin{pro} \label{5}
Let $(X, \{s_i\}_{i=1}^\infty, $ $\{ p_j
 \}_{j=1}^\infty)$ be  an infinite conformal uniformaly Dini IFS. Then the Ruelle operator $T$ is irreducible and
$$
{\dim} \{ h \in C(\overline J): Th = \varrho h, \ h \geq 0 \} \leq
1;
$$
if $h \geq 0$ is a $\varrho$-eigenfunction of $T$, then $h >0$.
\end{pro}

\medskip

\noindent{\bf Proof.}  For any given $f \in C(\overline J)$ with $f
\geq 0$ and $f\not\equiv 0$, let $V = \{ x \in \overline J: f(x) > 0
\}$.  For any $x \in \overline J$, by Proposition \ref{3.1}, there
exists a ${\mathbf u}_{0}$ such that $s_{{\mathbf u}_{0}} (x) \in
V$. Let $n_{0} = |{\mathbf u}_{0}|$, then
$$
T^{n_{0}}f(x) = \sum\limits_{|{\mathbf u}|=n_{0}}p_{s_{{\mathbf
u}}}(x)f(s_{{\mathbf u}}(x)) \geq p_{s_{{\mathbf
u}_{0}}}(x)f(s_{{\mathbf u}_{0}}(x)) > 0.
$$
This implies that  $T$ is irreducible.

For the dimension of the eigensubspace, we suppose that there exist
two independent  strictly positive $\varrho$-eigenfunctions $h_1, h_2
\in C(\overline J)$. Without loss of generality we assume that $0< h_1 \leq
h_2$ and $h_1(x_{0})=h_2(x_{0})$ for some $x_{0} \in \overline J$. Then $h =
h_2-h_1 (\geq 0)$ is a $\varrho$-eigenfunction of $T$ and  $h(x_{0}) =
0$. It follows that $T^{n}h(x_{0})=\varrho^{n} h(x_{0}) = 0,$
 which contradicts to the irreducibility of $T$. Hence
the dimension of the $\varrho$-eigensubspace is at most 1.

The strict positivity of $h$ follows directly from the
irreducibility of $T$.  \ \hfill $\Box$

\bigskip

\begin{pro} \label{P3.3}
Let $(X, \{s_i\}_{i=1}^\infty, $ $\{ p_j
 \}_{j=1}^\infty)$ be  an infinite conformal uniformaly Dini IFS
 and
$\varrho$ be the spectral radius of $T$. Then

{\em (i)} there exists  $0< h \in C(\overline J)$ such that
$Th=\varrho h;$

{\em (ii)} there exist $ A, B > 0$ such that
$$A \leq \varrho^{-n}T^{n}1(x) \leq B, \quad \forall n > 0.$$
\end{pro}

\bigskip

\noindent{\bf Proof.} (i) Let
\begin{eqnarray}\label{20}
\Phi(t)=\sum \limits_{j=0}^\infty \alpha (s^j t),
\end{eqnarray} where $\alpha(t)=\sup_{j \geq 1}\max_{|x-y|
\leq t}|\log p_j(x)-\log p_j(y)|$. From the uniformly Dini
continuity of $\{p_j\}_{j=1}^\infty$, it follows that $\Phi(\cdot)$
is well-defined and $\Phi(t)<+\infty.$

Let $C^+(\overline J):=\{f \in C(\overline J): f>0\}$, and set
$$F:=\{f \in C^+(\overline J): f(x) \leq f(y)e^{\Phi(|x-y|)}\}.$$
For any $f \in F$ and any $x, y \in \overline J$, we have
\begin{eqnarray} \label{P3.3-1}
&& Tf(x)=\sum \limits_{j=1}^\infty p_j(s_j(x)) f(s_j(x))\\
&&\leq\sum \limits_{j=1}^\infty p_j(s_j(y)) f(s_j(y))e^{\alpha (|x-y|)+\Phi(s|x-y|)} \nonumber\\
&&=Tf(y)e^{\alpha (|x-y|)+\Phi(s|x-y|)} \leq Tf(y)e^{\Phi(|x-y|)}.
\nonumber
\end{eqnarray}
We define $L: F \longrightarrow C(\overline J)$ by
$$L f(x):=\frac {Tf(x)}{\|Tf\|}.$$
By (\ref{P3.3-1}) we deduce
$$L f(x)\leq L f(y)e^{\Phi(|x-y|)}.$$
Let
$$F_0= F \bigcap \{ e^{-\Phi(1)}\leq f \leq 1 \}.$$
 It is easy to show that $F_0$ is a
convex compact subset of $ C(\overline J)$, and $L(F_0) \subset
F_0$. The Schauder fixed point  yields an $h \in F_0$ such that
$Lh=h$. This implies that $Th=\|Th\| h$. By the Proposition 3.2, we
have $\|Th\|=\varrho$. Hence $Th= \varrho h$.

(ii) The proof comes from Proposition \ref{P3.2} immediately. \hfill
$\Box$
\bigskip

We are now ready to prove our first result.\\

{\bf Proof of Theorem \ref{0.3}:} The proof is modified from the
paper \cite{LY}. We include the details here for the sake of
completeness. By Proposition \ref{P3.3}, there exists $0< h \in
C(\overline J)$ and constants $B \geq A>0$ such that $Th=\varrho h$
and
\begin{equation} \label{3.3-1}
A \leq \varrho^{-n}T^{n}1(x) \leq B, \qquad \forall \  n > 0.
\end{equation}

Let $C^+(\overline J):=\{f \in C(\overline J): f>0\}$, and let
$$D=\{ f \in C^+(\overline J):  \mbox{ there exists }\ c > 0  \ \mbox{such
that} \ f(x) \leq f(y) e^{c |x-y|}  \}. $$
 Then $D$ is dense in $C^{+}(\overline J)$.

Let
$$D_k:=\{f\in
C^+(\overline J):f(x)\leq f(y)e^{k\Phi(|x-y|)}\},$$ where $\Phi$ is
given in (\ref{20}). It is easy to see that $D \subseteq
\tilde{D}:=\bigcup _{k=1}^\infty D_k$,  hence $\tilde{D}$ is dense
in $C^+(\overline J)$.

For any $ f\in D_k$ and $x, y \in \overline J$, we have
\[ \begin{split}
T f(x) & =  \sum \limits_{j=1}^\infty p_j(s_j(x)) f(s_j(x))\\
& \leq  \sum \limits_{j=1}^\infty p_j(s_j(y))f(s_j(y))e^{k\alpha (|x-y|)+k\Phi(s|x-y|)}\\
& =  T f(y)e^{k\alpha (|x-y|)+ k \Phi(s|x-y|)} \\
& =  T f(y)e^{k\Phi(|x-y|)}.
\end{split} \]
It follows that $TD_k \subset D_k$.

For any $g \in C^+(\overline J)$, $f\in D_k$ and $n>0$, we have
\[ \begin{split}
 |\varrho ^{-n}T^n g(x) - \varrho ^{-n}T^n g(y)|\\
& \leq \|\varrho ^{-n}T^n f\||1- \frac {T^nf(y)} {T^nf(x)}|+2\|\varrho ^{-n}T^n\| \|f-g\|\\
& \leq  B \|f\|(e^{k\Phi(|x-y|)}-1)+2\|f-g\|,
\end{split} \] for any $x, y \in
\overline J$. By the assumptions on $D$ and $\Phi$, we can deduce
that for any $g \in C^+(\overline J)$, $\{\varrho ^{-n}T^n
g\}_{n=1}^\infty$ is a bounded equicontinuous sequence.

For any $f \in C(\overline J)$, we can choose $a>0$ such that
$f+a>0$. Then from the above prove process, it follows that the
sequences $\{\varrho ^{-n}T^n (f+a)\}_{n=1}^\infty$ and $\{\varrho
^{-n}T^n a\}_{n=1}^\infty$ are bounded equicontinuous, hence the
sequence $\{\varrho ^{-n}T^n f\}_{n=1}^\infty$ is also bounded
equicontinuous. We let
$$q_j(x)=\frac {p_j(s_j(x))h(s_j(x))}{\varrho h(x)}$$ and define an operator $P: C(\overline J)\rightarrow C(\overline
J)$ by
$$Pf(x)=\sum \limits_{j=1}^\infty q_j(x)f(s_j(x)).$$

For any $f \in C(\overline J)$, we have $T^nf=\varrho ^{n}h P^n
(f\cdot h^{-1})$. This implies that $\{P^n f\}_{n=1}^\infty$ is a
bounded equicontinuous sequence in $ C(\overline J)$. We know from
the Arzel\`{a}-Ascoli theorem that there exists $\tilde{f}\in
C(\overline J)$ and a subsequence $\{P^{n_i} f \}_{i=1}^\infty$ such
that $\lim\limits_{i\rightarrow \infty}\| P^{n_i} f-\tilde{f}\|=0$.

We claim that $\tilde{f}$ is a constant function and $\lim
\limits_{n\rightarrow \infty}\| P^{n} f-\tilde{f}\|=0$. For this we
let $\tau (g)=\min\limits_{x\in \overline J} g(x)$. Since $\sum
_{j=1}^\infty q_j(x)=1$, it is easy to see that $\tau
(\tilde{f})\leq \tau(P \tilde{f})$ and
\begin{equation} \label{4}
\tau (f)\leq \tau (P f)\leq \cdots \leq \tau (\tilde{f}).
\end{equation}
By taking the limit, we have $\tau (P\tilde{f})\leq \tau
(\tilde{f})$,  hence $\tau (P \tilde{f})= \tau (\tilde{f})$. For any
$n>0$, we select $x_n\in \overline J$ such that $P^n
\tilde{f}(x_n)=\tau (P^n \tilde{f})$. Then $\sum_{|{\mathbf
j}|=n}q_{s_{\mathbf {j}}}(x_n)=1$ implies that
$${\tilde{f}}(s_{\mathbf{j} }(x_{n})) = \tau({\tilde{f}}), \quad
\forall {\mathbf j}\in I^*, |{\mathbf j}| = n.$$
 Similarly there
exists $y_{n} \in \overline J$ such that
$${\tilde{f}}(s_{\mathbf
{j}}(y_{n})) = \eta({\tilde{f}}) := \max_{x\in \overline J}
{\tilde{f}}(x), \quad \forall {\mathbf j} \in I^*, |{\mathbf j}| =
n.$$
 We assume
${\mathbf j}_{n}=1 1 \cdots 1$, $|{\mathbf j}_{n}|=n$. Then
$$z:=\lim_{n
\rightarrow \infty}s_{{\mathbf j}_{n}}(x_{n})=\lim_{n \rightarrow
\infty}s_{{\mathbf j}_{n}}(y_{n}) \in \overline J.$$
Hence
$$\tau({\tilde{f}})=\lim_{n \rightarrow \infty }
{\tilde{f}}(s_{{\mathbf
j}_{n}}(x_{n}))={\tilde{f}}(z)=\lim_{n \rightarrow \infty}{\tilde{f}}(s_{{\mathbf
j}_{n}}(y_{n}))=\eta({\tilde{f}}).
$$
We duduce ${{\tilde{f}}}(x) \equiv \tau({{\tilde{f}}})$ is constant
function. By (\ref{4}) and the dual version for $\eta
({\tilde{f}})$, we have $\lim\limits_{n \rightarrow \infty} \|P^n
f-{\tilde f}\|=0$.

In particular, by taking $f = h^{-1}$, we see that $P^n (h^{-1})$
converges uniformly,  then $\varrho^{-n}T^{n}1$ converges uniformly.

To prove that
$$\lim_{n \to \infty} \| \varrho^{-n}T^{n}1 - h \|=0,$$
we let
$$
 f_{n}(x) = \frac{1}{n} \sum\limits_{i=0}^{n-1}\varrho^{-i}T^{i}1(x).
$$
From (\ref{3.3-1}), it follows that $\{ f_{n} \}_{n=1}^{\infty}$ is
bounded by $A$ and $B$ and is an equicontinuous subset of
$C(\overline J)$. By Arzel\`{a}-Ascoli theorem, we  assume without
loss of generality that there exists a $\tilde{h} \in C(\overline
J)$ such that $ \lim\limits_{n\rightarrow \infty}\|f_{n}-
\tilde{h}\|=0.$ We have
$$
\|T \tilde{h}-\varrho \tilde{h}\|=\lim\limits_{n\rightarrow \infty}
\|Tf_{n} - \varrho f_{n} \| \leq \lim\limits_{n\rightarrow
\infty}\frac{\varrho}{n} \| 1 - \varrho^{-n}T^{n}1 \| \leq
\lim\limits_{n\rightarrow \infty}\frac{\varrho}{n}\Big(1+B\Big)=0,
$$
i.e., $T \tilde{h} =\varrho \tilde{h}$ and also $\tilde{h} \geq A
>0$.   Proposition \ref{5}  implies that $\tilde{h} =c h$
for some $c>0$. Without loss of generality, we  assume that $h =
\tilde{h}.$ This implies that $\lim \limits_{n \rightarrow \infty}
\| \varrho^{-n}T^{n}1 - h \|=0.$ So $\lim \limits_{n\rightarrow
\infty}\| P^n (h^{-1})-1 \|=0$.

Now we define a function $\upsilon: C(\overline J) \rightarrow {\Bbb
R}$,
 $\langle \upsilon, f\rangle = \tau({\tilde{f}}) ( =\tilde{f} (x)$,
$x \in \overline J$). Then  $\upsilon $ is a bounded linear
functional on  $C(\overline J)$, and $\langle \upsilon, 1
\rangle=1$, $\langle \upsilon, h^{-1} \rangle=1$. From
 $$
 \langle \upsilon, Pf \rangle=\tau(P {\tilde{f}})=\tau({\tilde{f}}),
$$
we have $P^* \upsilon=\upsilon.$ Let $\mu: C(\overline J)
\rightarrow {\Bbb R}$ be defined by $\langle \mu, f \rangle =\langle
\upsilon, f h^{-1} \rangle$. Then $\langle \mu, 1 \rangle=\langle
\upsilon, h^{-1} \rangle=1$ and $\mu$ is a probability measure. It
is easy to see that $T^{*}\mu=\varrho \mu$ and $\langle \mu, h
\rangle=\langle \upsilon, 1 \rangle=1$. Hence for any  $f \in
C(\overline J)$, $\varrho^{-n}T^{n}f$ converges to  $\langle \mu, f
\rangle h$ in the supremum norm. Also it follows that for every $\xi
\in M(\overline J)$, $\varrho^{-n} T^{*n} \xi$ converges weakly to
$\langle \xi, h \rangle \mu$.

By Proposition \ref{5}, we conclude that the eigen-function $h$ is
unique. For the uniqueness of the eigen-measure, we note that if $\sigma \in M(\overline J)$ satisfies $T^{*} \sigma = \varrho
\sigma$ and $\langle \sigma, h \rangle = 1,$ then for every $f \in
C(\overline J)$,
$$\langle\sigma, f \rangle = \lim_{n\rightarrow
\infty} \ \langle \varrho^{-n} T^{*n} \sigma, f \rangle
=\lim_{n\rightarrow \infty}\ \langle \sigma, \varrho^{-n} T^{n}f
\rangle = \langle \sigma, \langle \mu, f \rangle h \rangle = \langle
\mu, f \rangle.$$
This implies that $\sigma = \mu$.  \hfill$\Box$


\bigskip
\section
{\bf SEPARATION PROPERTIES}\label{0.2}

 \bigskip
 \bigskip

 For any ${\mathbf i} \in I^*$, let
$$
J_{\mathbf i}=s_{\mathbf i}(J),\quad r_{\mathbf i}=\inf \limits_{x \in
U_0}|s'_{\mathbf i}(x)|, \quad R_{\mathbf i}=\sup \limits_{x \in
U_0}|s'_{\mathbf i} (x)|.
$$

\begin{lemma} \label {L4.1}  Let $X$ and $\{s_j\}_{j=1}^\infty$ be defined as
above. Then

{\em (i)} there exists  $c_1>1$ such that \begin{equation}
\label{4.1}
 R_{\mathbf j}\leq c_1 r_{\mathbf j}, \quad
\forall  {\mathbf j} \in I^*;
\end{equation}
\begin{equation} \label{4.2}
c_1^{-1}r_{\mathbf i}r_{\mathbf j} \leq r_{\mathbf ij} \leq c_1
r_{\mathbf i} r_{\mathbf j}, \quad \forall {\mathbf i}, {\mathbf j}
\in I^*;
\end{equation}

{\em (ii)} there exist $c_2 \geq c_1$ and $\delta >0$ such that for
any $x, y \in X$ with $|x-y| \leq \delta$,
\begin{equation} \label{4.3}
c_2^{-1}r_{\mathbf j} \leq \frac {|s_{\mathbf {j}}(x)-s_{\mathbf
{j}}(y)|}{|x-y|} \leq c_2 r_{\bf j}, \quad \forall {\mathbf {j}} \in
I^*;
\end{equation}

{\em (iii)} there exist $c_3 \geq c_2$ and $k_0$ such that for any
$x, y \in X$,
\begin{equation} \label{4.4}
|s_{\mathbf {j}}(x)-s_{\mathbf {j}}(y)| \leq c_3 r_{\mathbf
{j}}|x-y|, \quad \forall {\mathbf {j}} \in I^*, |{\mathbf j}|>k_0.
\end{equation}
\end{lemma}

\bigskip

\noindent{\bf Proof.}  (i) For each $j \in I$, let
$p_j(\cdot)=|s_j'(\cdot)|$ and still denote \eqref{20} by $\Phi(t)$.
For any $x, y \in U_0$, ${\mathbf j} \in I^*$, $|{\mathbf {j}}|=n$,
we have
\begin{eqnarray*}
&&|\log |s_{\mathbf {j}}'(y)|-\log|s_{\mathbf {j}}'(x)||\\
&& \leq \sum \limits_{i=1}^n \big|\log |s_{{j_i}}'(y_{i+1})|-\log|s_{
{j_i}}'(x_{i+1})| \big| \\
&&\leq \sum \limits_{i=1}^{n} \alpha (s^{n-i}|y-x|) \leq \Phi (1)<
\infty.
\end{eqnarray*}
where $y_{i}=s_{j_{i}} \cdots s_{j_n}(y)$, $y_{n+1}=y,$ $ y\in U_0$.
Consequently, we deduce (\ref{4.1}). And  the chain rule yields
(\ref{4.2}).

(ii) For any $x\in X$, there exist $\delta_x>0$ such that
$B(x,\delta_x)\subset U_0$. Since $X$ is locally connected, we can
assume that $B(x,\delta_x) \bigcap X$ is connected. Let $\delta$ be
the Lebesgue number. Then for any $x, y \in X$, if $|x-y| \leq
\delta$,  there exists $x'\in X$ such that $x, y \in
B(x',\delta_{x'})$. For such $x$ and $y$, we have $s_{\mathbf
{j}}(x),s_{\mathbf {j}}(y) \in B(y',\delta_{y'})$ for some $y' \in
X$. The self-similar property of $s_{\mathbf {j}}$ implies that
\begin{eqnarray*}
|s_{\mathbf {j}}(x)-s_{\mathbf {j}}(y)|\leq R_{\mathbf {j}}|x-y|
\leq c_1r_{\bf j}|x-y|.
\end{eqnarray*}

On the other hand, let  $u_{\mathbf j}(x):=s_{\mathbf {j}}^{-1}(x)$,
\  $\forall x \in B(y',\delta_{y'})\bigcap s_{\mathbf
{j}}(B(x',\delta_{x'}))$. Then
$$ R_{\bf j}^{-1}\leq |u_{\mathbf {j}}'(x)|\leq  r_{\mathbf {j}}^{-1}, \qquad
\forall x \in B(y',\delta_{y'}) \bigcap s_{\mathbf {j}}(B(x',
\delta_{x'})).$$
 Since $B(y', \delta_{y'}) \bigcap s_{\mathbf
{j}}(B(x', \delta_{x'}))$ is convex connected, similarly, by the
mean value theorem, we have
\begin{eqnarray*}
|u_{\mathbf {j}}(s_{\mathbf {j}}(x))-u_{\mathbf {j}}(s_{\mathbf
{j}}(y))|\leq r_{\mathbf {j}}^{-1}|s_{\mathbf {j}}(x)-s_{\mathbf
{j}}(y)|.
\end{eqnarray*}
Consequently, we have
$$r_{\mathbf {j}}|x-y|\leq |s_{\mathbf
{j}}(x)-s_{\mathbf {j}}(y)|\leq c_1r_{\mathbf {j}}|x-y|.$$

(iii) By the contractiveness of the infinite conformal IFS, we can
select integer $k_0$ such that
$$|s_{\mathbf {j}}(x)-s_{\mathbf {j}}(y)|\leq \delta, \quad \forall |{\mathbf j}|>k_0.
$$
The choice of the $\delta$ (the Lebesgue number) and the
self-similar property of $s_{\mathbf {j}}$ implies that
\begin{eqnarray*}
|s_{\mathbf {j}}(x)-s_{\mathbf {j}}(y)| \leq R_{\mathbf {j}}|x-y|
\leq c_3r_{\mathbf {j}}|x-y|.
\end{eqnarray*}
\hfill $\Box$
\bigskip

Let $\kappa$ be as given by (\ref{eq2.1}). We take $0<\varepsilon <
c_2^{-1} \cdot \min \{ \kappa, \delta \}$. By the assumption on $X$,
we have
\begin{equation} \label{4.5}
B(J,c_2\varepsilon)\subset X.
\end{equation}

For ${\mathbf j} \in I^*$, let $G_{\mathbf {j}}=s_{\mathbf
{j}}(B(J,\varepsilon))$. By (\ref{4.3}) and (\ref{4.5}), we have for
any $x\in J$,
\begin{equation} \label{4.6}
B(s_{\mathbf {j}}(x),c_2^{-1}\varepsilon r_{\mathbf {j}}) \subset
s_{\mathbf {j}}(B(x,\varepsilon)) \subset B(s_{\mathbf
{j}}(x),c_2\varepsilon r_{\mathbf {j}}).
\end{equation}
It follows that
\begin{eqnarray} \label{4.7}
&& B(J_{\mathbf j},c_2^{-1}\varepsilon r_{\mathbf j}) = \bigcup
\limits_{x\in J}B(s_{\mathbf j}(x),c_2^{-1}\varepsilon r_{\mathbf j
} )\subset G_{\mathbf j} =s_{\mathbf j}(\bigcup \limits_{x\in
J}B(x,\varepsilon ))\\
& &\subset  \bigcup \limits_{x\in J}B(s_{\mathbf
j}(x),c_2\varepsilon r_{\mathbf j}) =B(J_{\mathbf j},c_2\varepsilon
r_{\mathbf j}).\nonumber
\end{eqnarray}

For any two compact subsets $E, F$ of $R^d$, we define
\[ \begin{split}
&|E|  =  \sup\{|x-y|: x, y \in E \};\\
&D(E,F)  =  \inf \{|x-y|: x \in E, y \in F\};\\
&d(E,F) = \inf \{\varepsilon: E\subset B(F,\varepsilon ), \ F
\subset B(E,\varepsilon)\}.
 \end{split}\]

We say ${\mathbf {u}, \mathbf {v}} \in I^*$ are comparable if there
exists ${\mathbf {w}} \in I^*$ such that  $\mathbf {u} = \mathbf {v}
\mathbf {w}$
 or $\mathbf {v}=\mathbf {u} \mathbf {w}$.

Denote $F_n=\{1,2,\cdots ,n\}$. Let $F_n^*=\bigcup_{k \geq 1}
F_n^k$, and let
$$F_n^k=\{{\mathbf {u}}=u_1 u_2 \cdots u_k: u_i \in F_n,
1 \leq i \leq k \}.$$
 For any $n \in I$ and $0< r\leq 1$, we let
$$Q^n(r)=\{{\mathbf {v}}=v_1 v_2 \cdots v_m \in F_n^*: s_{\mathbf {v}}<r \leq
s_{v_1 v_2 \cdots v_{m-1}}\}.$$ Let $k_0$ be as given by Lemma
\ref{L4.1}(iii). For any ${\mathbf {w}}$ with $|{\mathbf {w}}|=k_0
+1$, we define
$$I^n({\mathbf {w}})=\{{\mathbf {v}} \in Q^n(|G_{\mathbf w}|): J_{\mathbf
{v}} \bigcap G_{\mathbf {w}} \neq \emptyset\}.$$
  Suppose $I^n(\mathbf
{w})$ is defined, and then for any $1 \leq j \leq n$, we define
$$I^n(j {\mathbf {w}})=M \bigcup N,$$
where $M=\{j{\mathbf {v}} :{\mathbf v}\in I^n({\mathbf {w}}) \}$ and
$$N=\{{\mathbf {v}} \in Q^n(|G_{j{\mathbf {w}}}|): v_1 \neq j,
J_{\mathbf {v}} \bigcap G_{j{\mathbf {w}}} \neq \emptyset\}.$$

\bigskip

\begin{lemma} \label{L4.2}
 For any fixed $n$, there exist constants $k_n$ and $c_{4}^{(n)}>0$
 such that
 $$\frac{1}{c_{4}^{(n)}} \leq \frac{r_{\mathbf {u}}} {r_{\mathbf
{v}}} \leq c_{4}^{(n)}, \quad \forall {\mathbf {u}}\in I^n({\mathbf
{v}}), \ |{\mathbf {v}}|\geq k_n.$$
\end{lemma}

\bigskip
\noindent{\bf Proof.} Let $k_0$ be as given in Lemma
\ref{L4.1}(iii). For any integer $n$, we select integer $k_n \geq
k_0$ such that
$$\min\{ |{\mathbf {u}}|: {\mathbf {u}} \in  I^n({\mathbf {v}}) \ \mbox{ and } \ |{\mathbf {v}}| \geq k_n \} > k_0.$$
We consider the following two case:

 (1) If $v_1 \neq u_1$, by the construction of
$N$, we have ${\mathbf {u}} \in Q^n(|G_{\mathbf {v}}|)$. Then
$$
r_{\mathbf {u}}< |G_{\mathbf {v}}|\leq r_{u_1 u_2 \cdots u_{m-1}}
\leq c_1 (r^{(n)})^{-1}r_{\mathbf {u}}, $$
 where
$r^{(n)}=\min_{1 \leq j \leq n} \{ r_j \}$. It follows from
(\ref{4.3}) that
$$c_2^{-1}\varepsilon r_{\mathbf {v}} \leq
|G_{{\mathbf {v}}}|.$$  Hence
$$
c_2^{-1}\varepsilon r_{\mathbf {v}} \leq |G_{{\mathbf {v}}}| \leq
c_1 (r^{(n)})^{-1}r_{\mathbf {u}}.
$$
Note that
$$r_{\mathbf {u}}<|G_{\mathbf {v}}|\leq c_2\varepsilon r_{\mathbf
{v}}+|J_{\mathbf {v}}| \leq c_3(2\varepsilon+|J|)r_{\mathbf {v}}.$$
From the  arguments above, we conclude that there exists $a>0$ such
that
$$\frac{1}{a} \leq \frac{r_{\mathbf u}}{r_{\mathbf v}}\leq a.$$

(2) If $v_1= u_1$, we write
$$
{\mathbf {v}}=v_1 v_2 \cdots v_l v_{l+1} \cdots v_n:=v_1 v_2 \cdots
v_l\mathbf {v}',$$
$${\mathbf {u}} = v_1 v_2 \cdots v_l u_{l+1} \cdots
u_n:=v_1 v_2 \cdots v_l\mathbf {u}'.
$$
 where $v_{l+1}\neq u_{l+1}$. From the construction of $M$,
  it follows that ${\mathbf {u}}' \in I^n({\mathbf {v}}')$. Similarly to the case of $(1)$, we can deduce that
$$a^{-1} \leq \frac{r_{\mathbf
{u}}'} {r_{\mathbf v}'}\leq a.$$ Together with (\ref{4.2}), this
implies that
$$(ac_1^2)^{-1}\leq \frac{r_{\mathbf {u}}} {r_{\mathbf {v}}} \leq
ac_1^2.$$ Let $c_{4}^{(n)}=ac_1^2$. Then the result follows.
\hfill$\Box$

\bigskip

\begin{pro} \label{4.3}
Let $\{s_j\}_{j=1}^\infty$ be an infinite conformal IFS.

{\em (i)} The system $\{s_j\}_{j=1}^\infty$ satisfies the finite
strong condition if it satisfies the open set condition;

 {\em (ii)}  $\sum \limits_{i=1}^\infty
R_i^d\leq c_1^d$ if the system $\{s_j \}_{j=1}^\infty$ satisfies
finite open set condition.
\end{pro}

\bigskip

\noindent{\bf Proof.} (i) For any $n$, from the assumption that the
 system $\{s_j \}_{j=1}^n$ satisfies open set condition,
 the result of Schief \cite{S} implies that the
system $\{s_j\}_{j=1}^n$ satisfies finite strong open set condition.
Hence there exists a nonempty bounded open set $U_n$ such that
$s_i(U_n) \subset U_n, \forall$ $ i \leq n$, and
$$s_i(U_n)\cap s_j(U_n)=\emptyset, \quad \forall \ i, j \leq n, i \neq
j.$$ Moreover $U_n \bigcap J_{F_n}\neq \emptyset$, where $J_{F_n}$
is the invariant set of the system $\{s_i \}_{i=1}^n$. It is easy to
see that $J_{F_n}\subset J$. Hence $U_n \bigcap J \neq \emptyset$.

(ii) Let $\lambda$ denote the Lebesgue measure on ${\Bbb{R}}^d$. For
any $n$, by definition of the finite open set condition, there
exists a nonempty bounded open set $U_n$ such that $s_i(U_n)\subset
U_n, \forall$ $ i \leq n$, and
$$s_i(U_n) \cap s_j(U_n)=\emptyset, \quad \forall\
i, j \leq n, i \neq j.$$
It follows that
$$\lambda (U_n) \geq \sum
\limits_{i=1}^n\lambda(s_i(U_n))=\sum \limits_{i=1}^n\int
\limits_{U_n}|s'_i(x)|^d d\lambda(x) \geq c_1^{-d}\sum
\limits_{i=1}^n R_i^d\lambda (U_n).$$ Since $\lambda (U_n)>0$, we
have
$$\sum
\limits_{i=1}^n R_i^d\leq c_1^d.$$ Hence $\sum \limits_{i=1}^\infty
R_i^d\leq c_1^d$. \hfill$\Box$

\bigskip

For any $t \geq 0$, let $\psi (t)=\sum \limits_{i\in I}R_i^t$, and
let $\theta =\inf \{t: \psi (t)<\infty \}$. For any $s >\theta $, we
define the {\it Ruelle operator} $T_s:C(\overline J)\rightarrow
C(\overline J)$ by
$$T_s f(x)=\sum \limits_{j=1}^\infty |s_j'(x)|^s \cdot f(s_j(x)).$$
Let $\varrho(T_s)$ denote the spectral radius of $T_s$.

\bigskip

\begin{pro} \label{New}
Assume $\psi(\theta)=\lim\limits_{t\rightarrow \theta}\psi(t)>1$ and
$\psi(\infty)=\lim\limits_{t\rightarrow \infty}\psi(t)<1$,  then
there exists a unique $a>\theta$ such that $\varrho(T_a) =1$.
\end{pro}

\bigskip

\noindent {\bf Proof.} We claim that the function
$$\varrho(T_a):=\lim_{n\rightarrow\infty}\|T^n_s1\|^{\frac 1
n}$$
 is continuous and strictly decreasing on $(\theta,+\infty)$. Indeed, it is easy to see that
$$T_s^n1(x)=\sum \limits_{|{\bf j}|=n}|s'_{\bf j}(x)|, \quad \forall x \in
\overline J.$$
This implies that
$$\|T_s^n1\|=\sup \limits_{x\in \overline J}\sum
\limits_{|{\bf j}|=n}|s'_{\bf j}(x)|^s.$$
 By Lemma 4.1, we have
$$\sum_{|{\bf j}|=n}r_{\bf j}^s\leq \sum \limits_{|{\bf j}|=n}|s'_{\bf j}(x)
|^s\leq \sum \limits_{|{\bf j}|=n}R_{\bf j}^s\leq c_1 \sum_{|{\bf
j}|=n}r_{\bf j}^s.$$ It follows that
$$\lim_{n\rightarrow\infty}\|T^n_s1\|^{\frac
1 n}=\lim_{n\rightarrow\infty}(\sum\limits_{|{\bf j}|=n}r_{\bf
j}^s)^{\frac 1 n}.$$

Let
$$f_n(s):=\frac 1 n \log
\sum\limits_{ |{\bf j}|=n}r_{\bf j}^s.$$
  For any $0\leq\lambda \leq 1$, $s_1,s_2\in (\theta,+\infty)$, by the
H$\ddot{o}$lder's inequality, we have
\[ \begin{split}
f_n(\lambda s_1+(1-\lambda){s_2}) & =  \frac 1 n \log \sum
\limits_{|{\bf j}|=n}r_{\bf j}^{\lambda {s_1}+(1-\lambda){s_2}}\\
& = \frac 1 n \log \sum \limits_{|{\bf j}|=n}r_{\bf j}^{\lambda
{}s_1}r_{\bf j}^{(1-\lambda) {s_2}} \\
& \leq  \frac 1 n \log(\sum \limits_{ |{\bf j}|=n}r_{\bf
j}^{s_1})^\lambda(\sum \limits_{ |{\bf
j}|=n}r_{\bf j}^{s_2})^{(1-\lambda)}\\
& = \lambda f_n( s_1) + (1-\lambda) f_n({s_2}).
\end{split} \]
We conclude that $f_n( \cdot )$ is convex. This implies that
$\lim\limits_{n\rightarrow\infty}\frac 1 n \log \sum\limits_{ |{\bf
j}|=n}r_{\bf j}^s$ is convex. Since the convex function on
$(\theta,+\infty)$ is continuous, we confirm that $\varrho(T_s)$ is
continuous.

If $s\in (\theta,+\infty)$, then there exists $i_0\in I$ such that
$r_{i_0}^s\geq r_j^s$ for every $j\in I$. We have
\begin{equation}\label{4.9}\sum_{|{\bf j}|=n}r_{\bf j}^{s+t}\leq \sum_{|{\bf
j}|=n}r_{\bf j}^sr_{\bf j}^t\leq \sum_{|{\bf j}|=n}r_{\bf j}^s c_1
r_{i_0}^{nt},\end{equation} for $t>0$. We can  check that
$\varrho(T_s)$ is strictly decreasing by using \eqref{4.9}.
Combining the assumption, we see that there exists a unique
$a>\theta$ such that $\varrho(T_a) =1$. \hfill $ \Box$

\bigskip

In the following we always let $a$ be the constant such that
$\varrho(T_a) =1$. From \cite{MU} we know that $a$ is the Hausdorff
dimension of $J$.

\bigskip

\begin{ex} Let $X=[\frac 1 4, \frac 1 2]$ and let $\{s_j \}_{j=1}^\infty=\{s_j:X\rightarrow X:j\geq
1\}$ be an infinite conformal IFS consisting of similarities
$s_j(x)=2^{-2j}x+2^{-j}-2^{-2j}$. Thus $|s'_j(\cdot)|=2^{-2j}$ and
the system satisfies the OSC. It is easy to see that $J$ is compact
and $a=\frac 1 2$. So from  Theorem \ref{0.3} we know that there
exist $h\in C(J)$  and  $\mu \in M(J)$ such that
\begin{eqnarray*}Th= h, \quad T^*\mu = \mu, \quad <\mu,
h>=1.
\end{eqnarray*}
\end{ex}

\bigskip

\begin{pro} \label{P4.4}
 $H^a ( J)\leq H^a ( \overline J)<\infty$.
\end{pro}

\bigskip

 \noindent{\bf Proof.} For any fixed $z \in \overline J$,
by (\ref{4.4}), there exist $k_0$ and $|{\mathbf {j}}|>k_0$ such that
$|J_{\mathbf {j}}| \leq c |s'_{\mathbf {j}}(z)|,$ so we have
\begin{eqnarray*}
 \sum \limits_{{\mathbf {j}}=n}|J_{\mathbf {j}}|^a\leq c^a\sum \limits_{{\mathbf {j}}=n}|s'_{\mathbf {j}}(z)|^a.
\end{eqnarray*}
From Theorem \ref{0.3}, we know that there exists a unique $h(z)\in
C(\overline J)$ such that
$$\lim \limits_{n\rightarrow \infty}\sum \limits_{{\mathbf
{j}}=n}|s_{\mathbf {j}}'(z)|^a=\lim \limits_{n\rightarrow \infty}
T_a^n 1(z)=h(z).$$
 This implies $H^a ( J) \leq H^a ( \overline
J)<\infty$.\hfill$\Box$

\bigskip

 \begin{lemma} \label{L4.5} \cite{FL} Let $a>0$ be a constant. Let $w$ be conformal
and invertible and  $D$ be a Borel subset in the domain of $w$ and
$0<H^a (D)<\infty$. Then we have the following change of variable
formula:
$$H^a (w(D))=\int \limits_D |w'(x)|^adH^a (x).$$
\end{lemma}

\bigskip

 \begin{lemma} \label{L4.6}
 Suppose $0<H^a (J)<\infty$, then
$$
H^a (J_{\bf u}\bigcap J_{\bf v})=0 \ \mbox{ for any incomparable } \
{\mathbf {u}}, {\mathbf {v}} \in I^*.
$$
\end{lemma}

\bigskip

\noindent{\bf Proof.} Note that $\{| s_j'(\cdot)| \}_{j=1}^\infty$
is uniformly Dini continuous. Since $T_a$ has spectral radius 1, by
Theorem \ref{0.3}, there exists $0<h \in C(\overline J)$ such that
$$h(x)=\sum \limits_{j=1}^\infty |s_j'(x)|^a h(s_j(x)), \quad \forall x \in \overline J.$$
Let $\tilde{h} :=h|_J$. Then $0 < \tilde{h} \in C( J)$ and
$\tilde{h}(x)=\sum \limits_{j=1}^\infty |s_j'(x)|^a
\tilde{h}(s_j(x))$. Hence
\begin{eqnarray*} &&\sum
\limits_{j=1}^\infty \int \limits_{J_j} \tilde{h}(x)d{H^a }(x) \geq
\int \limits_J \tilde{h}(x) d{H^a}(x) \\
&&=\int \limits_J\sum
\limits_{j=1}^\infty|s_j'(x)|^a \tilde{h}(s_j(x)) d{H^a}(x)\\
&&=\sum\limits_{j=1}^\infty\int \limits_{J_j} \tilde{h}(x)
d{H^a}(x).
\end{eqnarray*}
(The last equality follows  Lemma \ref{L4.5}.)
 This  implies that
$H^a(J_i\bigcap J_j)=0$ for any $i\neq j$. It follows immediately
that $H^a (J_{\bf u}\bigcap J_{\bf v})=0$ for any incomparable $
{\bf u},{\bf v}\in I^*.$\hfill$\Box$

\bigskip
{\bf Proof of  Theorem  \ref{0.7}:} Let $c_1,$ $c_2$, $c_3$ and
$\delta$ be as given in Lemma \ref{L4.1}. And let $\iota>0$ satisfy
the condition: $2^{-1}c_1^{-a}>\iota$. There exists an open covering
$V_1, \cdots, V_n$ of $\overline J$ such that
\begin{equation} \label{4.8}
V:=\bigcup \limits_{i=1}^n V_i \supset \overline J, \quad
\delta':=D(\overline J,V^c)<\delta
\end{equation}
and
$$
H^a(V)\leq\sum \limits_{i=1}^n |V_i|^a\leq(1+\iota)H^a( \overline
J)\leq (1+\iota)H^a (J).
$$

For any given  ${\mathbf {u}}, {\mathbf {v}} \in I^n(\mathbf {w})$,
we let $\overline J_{\mathbf {u}}=\overline {s_{\mathbf {u}}(J)}$
and $\overline J_{\mathbf {v}}=\overline {s_{\mathbf {v}}(J)}$. We
assume that $H^a (J_{\mathbf u}) \leq H^a (J_{\mathbf v})$. Then for
any given $\varepsilon>0$ satisfying $c_1^a \iota<\varepsilon<1$, we
have $\varepsilon H^a (J_{\mathbf {u}}) \leq H^a (J_{\mathbf {v}})$.
We claim that $d(\overline J_{\mathbf {u}}, \overline J_{\mathbf
{v}} ) \geq c_2^{-1}\delta' r_{\mathbf {u}}$. Otherwise $d(\overline
J_{\mathbf {u}}, \overline J_{\mathbf {v}})<c_2^{-1}\delta'
r_{\mathbf {u}}$, by (\ref{4.3}) and (\ref{4.8}), we have
$$D(\overline J_{\mathbf {u}}, s_{\mathbf
{u}}(V)^c) \geq c_2^{-1} \delta' r_{\mathbf {u}}.$$
This implies
that $J_{\mathbf v}\subset \overline J_{\mathbf v} \subset
s_{\mathbf u} (V)$. By Lemma \ref{L4.6}, we have
$$(1+\varepsilon)H^a (J_{\mathbf
u}) <H^a (J_{\mathbf u})+H^a (J_{\mathbf v})=H^a (J_{\mathbf
u}\bigcup J_{\mathbf v})\leq H^a (s_{\mathbf u}(V)).$$
Thus
\begin{eqnarray*}
&&\varepsilon r_{\mathbf u}^a H^a( J)\leq\varepsilon H^a( J_{\mathbf
u})< H^a( s_{\mathbf u}(V\backslash J))\\&&\leq {( c_1r_{\mathbf
u})}^aH^a(V\backslash J)
 <{( c_1r_{\mathbf u })}^a\iota H^a( J).
\end{eqnarray*}
This implies that $\varepsilon<c_1^a\iota$, and it contradicts to the choice of
$\varepsilon$. Hence, there exists $\delta_0>0$ such that for any ${\bf
{w}} \in I^*$,
\begin{equation} \label{(4.9)}
d(\overline J_{\bf u}, \overline J_{\mathbf {v}}) \geq \delta_0
r_{\mathbf {w}}, \quad  \forall \ {\mathbf {u}}, {\mathbf {v}} \in
I^n({\mathbf {w}}), {\mathbf {u}} \neq {\mathbf {v}}.
\end{equation}

Let $\delta_1=(3c_3c_4^{(n)})^{-1}\delta_0$. We can find a finite set
$Z \subset \overline J$ such that
$$\overline J \subset \bigcup
\limits_{x\in Z}B(x,\delta_1).$$ For any ${\mathbf {w}} \in I^*$
with $|{\mathbf {w}}| > k_n$, by (\ref{(4.9)}) there exists $x \in
\overline J$ such that
$$|s_{\mathbf {u}}(x)-s_{\mathbf {v}}(x)| \geq \delta_0r_{\mathbf {w}}, \quad \forall \ {\mathbf {u}}, {\mathbf {v}} \in I^n(\mathbf {w}), {\mathbf {u}} \neq {\mathbf {v}}.$$
For such $x$ there exists a $z$ such that $|x-z| < \delta_1$. By
(\ref{4.4}) and the choice of $k_n$ in Lemma \ref{L4.2}, we have
$$|s_{\mathbf {u}}(x)-s_{\mathbf {u}}(z)| \leq \frac 1 3 \delta_0 r_{\mathbf {w}} ,\quad |s_{\mathbf {v}}(x)-s_{\mathbf {v}}(z)|\leq \frac 1 3 \delta_0 r_{\mathbf {w} }.$$
This implies that
\begin{equation} \label{4.10}
|s_{\mathbf {u}} (z)-s_{\mathbf {v}} (z)| \geq \delta_0 r_{\mathbf
{w}}/3.
\end{equation}
 For each $z\in Z$, let
$$T_z^n=\{{\mathbf {u}} \in I^n({\mathbf {w}}): \exists {\mathbf {v}} \in I^n({\mathbf {w}}) \  \mbox{ such that } (\ref{4.10}) \ \mbox{ holds } \}.
$$
Together with (\ref{4.10}), it implies that
$$I^n({\mathbf {w}})=\bigcup \limits_{z \in Z} T_z^n.$$
 For each $z$, the sets
$$\{B(s_{\mathbf u} (z),\delta_0 r_{\mathbf w} /6):{\mathbf u}\in
T^n_z\}$$
 are disjoint and  contained in $B( G_{\mathbf {w}},|\overline
J_{\mathbf {u}}|+\delta_0 r_{\mathbf {w}} /6)$. By (\ref{4.4}) and
Lemma \ref{L4.2}, there exist $x \in \overline J$ and $c>0$ such
that
$$B(G_{\mathbf {w}}, |\overline J_{\mathbf {u}}|+\delta_0 r_{\mathbf
{w} }/6) \subset B(x, cr_{\mathbf {w}}).$$ We deduce that there
exists $k$ such that $\max\limits_{z \in Z} \sharp T_z^n \leq k$,
then we have
$$\sharp I^n({\mathbf
{w}}) \leq \sharp Z \cdot \max_{z \in Z} \sharp T_z^n \leq \sharp Z
\cdot k.$$ If ${\mathbf {w}} \in I^*$ with $|{\mathbf {w}}|>k_n$ and
$$\sharp I^n({\mathbf {w}})=\sup_{|{\mathbf {w}}| \geq k_n}
\sharp I^n({\mathbf {w}}),
$$
we prove that
\begin{equation}\label{21}
I^n({\mathbf {i}} {\mathbf {w}})=\{{\mathbf {i}} {\mathbf {j}}:
{\mathbf {j}} \in I^n({\mathbf {w}})\} \quad \forall \ {\mathbf {i}}
\in F_n^*.
\end{equation}
 By the maximality of $\mathbf w$, we only prove the following:
$$\{ {\mathbf {i}} {\mathbf {l}}: {\mathbf {l}} \in I^n({\mathbf
{w}})\} \subset I^n({{\mathbf {i}} {\mathbf {w}}}).$$
The definition of
$I^n({\mathbf {w}})$ implies that
$$I^n(j {\mathbf {w}}) \supset \{j {\mathbf {v}}: {\mathbf {v}} \in I^n({\mathbf {w}})\}, \quad j=1, 2, \cdots, n.$$
We conclude that
$$\{{\mathbf {i}} {\mathbf {l}}: {\mathbf {l}} \in I^n({\mathbf
{w}})\} \subset I^n({\mathbf {i}} {\mathbf {w}}).$$
 For any fixed $1 \leq l \leq n$,  ${\mathbf {v}}=v_1 \cdots v_n \in I^*$, $v_1 \neq l$,
 we consider the family
$${\mathfrak{J}}_l=\{J_{\mathbf l}:{\mathbf l}\in Q^n(|{\bf vw}|),\
l_1=l\}.
$$
Then ${\mathfrak{J}}_l$ is a cover of $J_ l$. Since $v_1\neq l$, then
${\mathbf {l}} \notin I^n({\mathbf {vw}})$. By the construction of $N$, we have
$J_{\mathbf {l}} \bigcap G_{\mathbf {vw}}=\emptyset$. Hence by (\ref{4.7}) we
have
$$D(J_{\mathbf {l}}, J_{\mathbf {vw}}) \geq c_2^{-1}\varepsilon r_{\mathbf
{vw}},$$
which implies
\begin{equation} \label {4.11}
D(J_l,J_{\bf {vw}})\geq c_2^{-1}\varepsilon r_{\bf {vw}},  \quad
l\neq v_1.
\end{equation}

Now we let $U_n=\bigcup \limits _{{\mathbf u}\in F_n^*} G_{\bf
{uw}}^*$, where  $G_{\mathbf u}^*=s_{\mathbf {u}} \big(B(J,
2^{-1}c_2^{-2}\varepsilon ) \big)$. We claim that the set $U_n$
satisfies the condition of finite SOSC. Indeed, $U_n$ is an open set
and for each $i \in F_n$,
$$s_i(U_n)=\bigcup \limits_{{\mathbf {u} }\in F_n^*} s_i(G_{{\mathbf
{u}} {\mathbf {w}}}^*)=\bigcup \limits _{{\mathbf {u}} \in
F_n^*} G_{i{\mathbf {u}}{ \mathbf {w}}}^* \subset U_n.
$$
We clain that for each $ i, j \in F_n$, $i\neq j$, $s_i(U_n)\bigcap
s_j(U_n)=\emptyset$. Otherwise, there exist ${\mathbf u}, { \mathbf
v }\in F_n^*$ such that  $G_{i{\mathbf u}{ \mathbf w}}^*\bigcap
G_{j{\bf{vw}}}^*\neq \emptyset$. We assume that
$$r_{i{\mathbf u}{\mathbf
w}}\geq r_{j{\bf{vw}}}.$$
 Let $y$ be in the intersection, then there
exist $y_1 \in J_{i{\mathbf u }{\mathbf w}}$ and $y_2 \in
J_{j{\bf{vw}}}$ such that
$$
d(y,y_1)<c_2\cdot \frac 1  {2c_2^2}\varepsilon\cdot r_{i{\mathbf u}{ \mathbf w}
} \leq \frac {c_2^{-1} \varepsilon } 2 r_{i{\mathbf {u}}{ \mathbf {w}},}$$
$$d(y, y_2) < c_2 \cdot \frac 1 {2c_2^2} \varepsilon \cdot r_{j{\mathbf {u}}
{\mathbf {w}}} \leq \frac {c_2^{-1} \varepsilon } 2 r_{j{\mathbf {u}}{
\mathbf w}}.
$$
Hence
$$D( J_{i{\mathbf {uw}}}, J_{j}) <c_2^{-1} \varepsilon r_{i{\mathbf
{uw}}},
$$
which contradicts to (\ref{4.11}). And it is easy to see that $U_n
\bigcap J \neq \emptyset$. This completes the proof. \hfill$\Box$

\bigskip
\bigskip

\end{document}